\documentclass[10pt]{amsart}
\usepackage{amsfonts,amssymb,amscd,amsmath,enumerate,verbatim,calc}
\usepackage{color,soul}
\usepackage{mathtools}
\usepackage{graphicx}
\everymath{\displaystyle}
\usepackage{fancyhdr}
\usepackage{tikz}
\usetikzlibrary{shapes.geometric, arrows}
\tikzstyle{arrow} = [thick,->,>=stealth]
\tikzstyle{process} = [rectangle, minimum width=4cm, minimum height=2cm, text centered, text width=2.5cm, draw=green, fill=green!10]
\theoremstyle{plain}

\theoremstyle{definition}

\begin{document}
\title{SC*-Normal spaces and some functions}
\author{Neeraj kumar Tomar}
\email{neer8393@gmail.com}
\address{Department of Applied Mathematics, Gautam Buddha University, Greater Noida, Uttar Pradesh 201312, India}

\author{Fahed~Zulfeqarr}
\email{fahed@gbu.ac.in}
\address{Department of Applied Mathematics, Gautam Buddha University, Greater Noida, Uttar Pradesh 201312, India}

\author{M.C Sharma}
\email{sharmamc2@gmail.com}
\address{Department of Mathematics , N.R.E.C College Khurja,203131, India}

\date{\today}
	
\begin{abstract}
In this paper, we introduce and explore a new class of topological spaces termed as \( SC^* \)-normal spaces, defined via \( SC^* \)-open sets. The concept of \( SC^* \)-normality is analyzed in relation to classical notions such as normal spaces and \( g \)-normal spaces. We further define and examine generalized forms of \( SC^* \)-closed functions, including \( SC^* \)-generalized closed mappings, and investigate their fundamental properties. Several characterizations of \( SC^* \)-normal spaces are established, and various preservation theorems under different types of functions are also presented.
\end{abstract}
\maketitle

\section*{Introduction} 
Normality is an important topological property and hence it is of significance both from intrinsic interest as well as from applications view point to obtain factorizations of normality in terms of weaker topological properties. In 1937, Stone\cite{stone1937applications} introduced the concept of regular-open sets. In 1963, Levine\cite{levine1963semi} introduced the notion of semi-open sets and obtained their properties. Levine\cite{levine1970generalized} initiated the investigation of $g$-closed sets in topological spaces, since then many modifications of $g$-closed sets were defined and investigated by a large number of topologists.

In this paper, we organize our study into five main sections.\\
\textbf{Section 1.} We begin by introducing the concepts of generalized $SC^*$-closed and $SC^*$-open sets, along with their variants: $gSC^*$-closed, $gSC^*$-open, $SC^*g$-closed, and $SC^*g$-open sets. Various examples are provided to illustrate these definitions, and we examine how these sets relate to other known types of open sets.\\
\textbf{Section 2.} We introduce the notion of $SC^*$-normal spaces and provide examples to illustrate this concept. Additionally, we establish an equivalent characterization, referred to as the fundamental theorem of $SC^*$-normal spaces. We also examine how $SC^*$-normality interacts with other topological properties and explore conditions under which such spaces behave analogously to classical normal spaces.\\
\textbf{Section 3.} This section focuses on the study of strongly $SC^*$-open, strongly $SC^*$-closed, and almost $SC^*$-irresolute functions. We develop several results, including theorems and lemmas, that connect these types of functions with $SC^*$-normal spaces. Their structural significance is analyzed, and we also investigate how these functions preserve or reflect topological features under certain mappings.\\
\textbf{Section 4.} Here, we explore various generalized functions such as $gSC^*$-closed, $SC^*g$-closed, quasi $SC^*$-closed, $SC^*$-$SC^*g$-closed, $SC^*$-$gSC^*$-closed, and almost $gSC^*$-closed functions. We also discuss some of their key properties and structural behavior. The relationships among these functions are examined through comparative analysis, and several examples are constructed to clarify their distinct characteristics.
\newpage
\textbf{Section 5.} In this final section, we present several preservation results and characterizations related to $SC^*$-normal spaces. The section includes proofs of relevant theorems, lemmas, and corollaries to support the theoretical framework. Furthermore, we explore how certain properties of $SC^*$-normal spaces are retained under specific function types, and highlight directions for further research in this area.

\section{Preliminaries and Notations} 
In this paper, we consider topological spaces $(X, \tau)$, $(Y, \sigma)$, and $(Z, \gamma)$, where no separation axioms are assumed unless stated otherwise. Functions between these spaces are denoted by $f: X \rightarrow Y$ and $g: Y \rightarrow Z$. For any subset $A \subseteq X$, the symbols $cl(A)$ and $int(A)$ represent the closure and interior of $A$, respectively.\\

We now define some basic notions which will be used throughout. For a good understanding of them, readers are referred to see \cite{levine1963semi, malathi2017pre, stone1937applications}.

\subsection{Definition:}
A subset $A$ of a topological space $X$ is said to be: \\
$(1)$ regular open\cite{stone1937applications} if $A=int(cl(A))$.\\ 
$(2)$ semi open\cite{levine1963semi} if $A\subset cl(int(A))$.\\ 
$(3)$ $c^*$-open\cite{malathi2017pre}  if $int(cl(A))\subset A\subset cl(int(A))$. 

\noindent A set is said to be regular-closed (respectively, semi-closed or $c^*$-closed) if it is the complement of a regular-open (respectively, semi-open or $c^*$-open) set. Given a subset $A$ of a topological space, the $c^*$-closure (respectively, semi-closure) of $A$, denoted by $c^*$-$cl(A)$ (respectively, $scl(A)$), is defined as the smallest $c^*$-closed (respectively, semi-closed) set containing $A$, i.e., the intersection of all such sets that include $A$. Similarly, the $c^*$-interior (respectively, semi-interior) of $A$, written as $c^*$-$int(A)$ (respectively, $s$-$int(A)$), is the largest $c^*$-open (respectively, semi-open) set contained within $A$, obtained by taking the union of all such subsets of $A$.

\subsection{Definition} A subset $A$ of a topological space $X$ is said to be $SC^*$-closed\cite{chandrakala2024closed} if $scl(A)\subset U$ whenever $A\subset U$ and $U$ is $c^*$-open in $X$. The complement of $SC^*$-closed set is said to be $SC^*$-open.

For any subset $A$ of a topological space, the $SC^*$-closure of $A$, denoted by $SC^*$-$cl(A)$, is defined as the smallest $SC^*$-closed set containing $A$, that is, the intersection of all $SC^*$-closed sets which include $A$. Similarly, the $SC^*$-interior of $A$, written as $SC^*$-$int(A)$, is the largest $SC^*$-open set contained in $A$, formed by taking the union of all such subsets.\\
The collection of all $SC^*$-open sets in a space $X$ is denoted by $SC^*O(X)$. Correspondingly, the notations $SC^*C(X)$, $RO(X)$, $RC(X)$, $SO(X)$, and $SC(X)$ are used to represent the families of $SC^*$-closed, regular open, regular closed, semi-open, and semi-closed sets in $X$, respectively.

\subsection{Definition} A subset $A$ of a topological space $(X,\tau )$ is said to be 
\begin{flushleft}
    
$(1)$ {$g$-closed}\cite{levine1970generalized}  if $cl(A)\subset U$ whenever $A\subset U$ and 
$U\in \tau$. 

$(2)$ {generalized $SC^*$-closed} (briefly $gSC^*$-closed)  if  $SC^*$-$cl(A)\subset U$ whenever $A\subset U$ and $U\in\tau$.

$(3)$ {$SC^*$generalized-closed} (briefly $SC^*g$-closed)  if  $SC^*$-$cl(A)\subset U$ whenever $A\subset U$ and $U\in SC^*O(X)$.\\

A set is called $g$-open (respectively, $gSC^*$-open or $SC^*g$-open) if its complement is $g$-closed (respectively, $gSC^*$-closed or $SC^*g$-closed).
\end{flushleft}

\subsection{Remark.} 
The relationships among various types of closed sets can be expressed through the following chain of implications:

\begin{flushleft}
    closed $\Rightarrow$ $SC^*$-closed  $\Leftrightarrow$ $SC^*g$-closed  $\Leftrightarrow$  $gSC^*$- closed \\$\Downarrow$\\ $g$-closed
\end{flushleft}

The relationships among various types of closed sets can be expressed through the following chain of implications:
\subsection{Example} Let $X=\{q,r,s,t\}$ and $\tau = \{\phi,\{r,t\},\{q,r,t\},\{r,s,t\},X\}$. Then

\begin{flushleft}
    $(1)$ closed sets in $(X,\tau)$ are $\phi$, $X$, $\{q\}$, $\{s\}$, $\{q,s\}$.

 $(2)$ $g$-closed set in $(X,\tau)$ are $\phi$, $\{q\}$, $\{s\}$, $\{q,s\}$, $\{q,r,s\}$, $\{q,s,t\}$.

 $(3)$ $SC^*$-closed set in $(X,\tau)$ are $\phi$, $X$, $\{q\}$, $\{r\}$, $\{s\}$, $\{t\}$, $\{q,r\}$, $\{q,s\}$, $\{q,t\}$, $\{r,s\}$, $\{r,t\}$, $\{s,t\}$, $\{q,r,s\}$, $\{q,r,t\}$, $\{q,s,t\}$, $\{r,s,t\}$.

$(4)$ $gSC^*$-closed set in $(X,\tau)$ are $\phi$, $X$, $\{q\}$, $\{r\}$, $\{s\}$, $\{t\}$, $\{q,r\}$, $\{q,s\}$, $\{q,t\}$, $\{r,s\}$, $\{r,t\}$, $\{s,t\}$, $\{q,r,s\}$, $\{q,r,t\}$, $\{q,s,t\}$, $\{r,s,t\}$.

$(5)$ $SC^*g$-closed set in $(X,\tau)$ are $\phi$, $X$, $\{q\}$, $\{r\}$, $\{s\}$, $\{t\}$, $\{q,r\}$, $\{q,s\}$, $\{q,t\}$, $\{r,s\}$, $\{r,t\}$, $\{s,t\}$, $\{q,r,s\}$, $\{q,r,t\}$, $\{q,s,t\}$, $\{r,s,t\}$.
\end{flushleft}

\subsection{Example} Let $X=\{q,r,s,t\}$ and $\tau = \{\phi,\{q\},\{r\},\{q,r\},\{q,s\},\{q,r,s\},X\}$. Then

\begin{flushleft}
    $(1)$ closed sets in $(X,\tau)$ are $\phi,\{t\},\{r,t\},\{s,t\},\{q,s,t\},\{r,s,t\},X$

$(2)$ $g$-closed set in $(X,\tau)$ are $\phi,X,\{t\},\{q,t\},\{r,t\},\{s,t\},\{q,r,t\},\{q,s,t\},\{r,s,t\}$.

$(3)$ $SC^*$-closed set in $(X,\tau)$ are $\phi$, $X$, $\{q\}$, $\{r\}$, $\{s\}$, $\{t\}$, $\{q,r\}$, $\{q,s\}$, $\{q,t\}$, $\{r,s\}$, $\{r,t\}$, $\{s,t\}$, $\{q,r,s\}$, $\{q,r,t\}$, $\{q,s,t\}$, $\{r,s,t\}$.

$(4)$ $gSC^*$-closed set in $(X,\tau)$ are $\phi$, $X$, $\{q\}$, $\{r\}$, $\{s\}$, $\{t\}$, $\{q,r\}$, $\{q,s\}$, $\{q,t\}$, $\{r,s\}$, $\{r,t\}$, $\{s,t\}$, $\{q,r,s\}$, $\{q,r,t\}$, $\{q,s,t\}$, $\{r,s,t\}$.

$(5)$ $SC^*g$-closed set in $(X,\tau)$ are $\phi$, $X$, $\{q\}$, $\{r\}$, $\{s\}$, $\{t\}$, $\{q,r\}$, $\{q,s\}$, $\{q,t\}$, $\{r,s\}$, $\{r,t\}$, $\{s,t\}$, $\{q,r,s\}$, $\{q,r,t\}$, $\{q,s,t\}$, $\{r,s,t\}$.
\end{flushleft}

 \section{$SC^*$-normal spaces}
\subsection{Definition} A topological space $X$ is said to be normal\cite{tietze1923beitrage} (resp. $g$- normal\cite{munsi1986sep}, $SC^*$-normal) if for any pair of disjoint closed sets $A$ and $B$, there exist open (resp. $g$-open, $SC^*$-open) sets $U$ and $V$ such that $A\subset U$ and $B\subset V$.

\subsection{Remark.}
In any topological space $X$, the following sequence of implications holds among different notions of normality:

\vspace{3.7mm}

\begin{center}
    normal $\Rightarrow$ $g$-normal $\Rightarrow$ $SC^*$-normal
\end{center}

\vspace{3.9mm}

However, these implications are not reversible in general. This fact will be supported by appropriate counterexamples presented in the following section.

\subsection{Example.}
Consider the finite set $X = \{q, r, s, t\}$ with the topology
\[
\tau = \{\emptyset, \{q\}, \{r\}, \{q, r\}, \{s, t\}, \{q, s, t\}, \{r, s, t\}, X\}.
\]
Let $A = \{q\}$ and $B = \emptyset$. Clearly, both $A$ $\&$ $B$ are closed in $X$ and are disjoint. We can find disjoint open sets $U = \{q, s, t\}$ $\&$ $V = \{r\}$ s.t., $A \subseteq U$ $\&$ $B \subseteq V$. 
This confirms that $X$ is a normal space. Since normality implies $g$-normality, and all open sets in this topology are also $SC^*$-open, it follows that $X$ is also $g$-normal $\&$ $SC^*$-normal.

\subsection{Theorem.}
For a topological space $X$, the following conditions are equivalent:

\begin{enumerate}
    \item The space $X$ is $SC^*$-normal.\\
    \item Given any two open sets $U$ $\&$ $V$ s.t. $U \cup V = X$, there exist $SC^*$-closed sets $A$ $\&$ $B$ with $A \subseteq U$, $B \subseteq V$, $\&$ $A \cup B = X$.\\
    \item For every closed subset $H$ of $X$ and every open set $K$ containing $H$, $\exists$ an $SC^*$-open set $U$ s.t.,
    \[
    H \subseteq U \subseteq SC^*\text{-}cl(U) \subseteq K.
    \]
\end{enumerate}
\begin{flushleft}
\textbf{Proof.}
\textbf{(1) $\Rightarrow$ (2):}  
Assume $X$ is an $SC^*$-normal space, and let $U$ $\&$ $V$ be open subsets of $X$ s.t., $U \cup V = X$. Then the complements $X \setminus U$ $\&$ $X \setminus V$ are disj., closed sets. By the definition of $SC^*$-normality, there exist disj., $SC^*$-open sets $U_1$ $\&$ $V_1$ satisfying $X \setminus U \subseteq U_1$ $\&$ $X \setminus V \subseteq V_1$. Define $A = X \setminus U_1$ $\&$ $B = X \setminus V_1$. It follows that $A$ $\&$ $B$ are $SC^*$-closed subsets s.t., $A \subseteq U$, $B \subseteq V$, $\&$ $A \cup B = X$.

\vspace{2mm}

\textbf{(2) $\Rightarrow$ (3):}  
Let $H$ be a closed subset of $X$, and let $K$ be an open set with $H \subseteq K$. Then $X \setminus H$ $\&$ $K$ are open sets whose union covers $X$. By hypothesis, $\exists$ $SC^*$-closed sets $M_1 \subseteq X \setminus H$ $\&$ $M_2 \subseteq K$ s.t., $M_1 \cup M_2 = X$. Taking complements, set $U = X \setminus M_1$ $\&$ $V = X \setminus M_2$. Then $U$ $\&$ $V$ are $SC^*$-open, disj., and satisfy:
\[
H \subseteq U,\quad X \setminus K \subseteq V,\quad \text{and} \quad U \subseteq X \setminus V \subseteq K.
\]
Since $X \setminus V$ is $SC^*$-closed, we conclude that $SC^*\text{-}cl(U) \subseteq X \setminus V \subseteq K$, which gives:
\[
H \subseteq U \subseteq SC^*\text{-}cl(U) \subseteq K.
\]

\vspace{2mm}

\textbf{(3) $\Rightarrow$ (1):}  
Let $H_1$ $\&$ $H_2$ be two disj., closed subsets of $X$. Define the open set $K = X \setminus H_2$, which contains $H_1$. By assumption, $\exists$ $SC^*$-open set $U$ s.t:
\[
H_1 \subseteq U \subseteq SC^*\text{-}cl(U) \subseteq K.
\]
Now, define $V = X \setminus SC^*\text{-}cl(U)$. Then $V$ is $SC^*$-open and contains $H_2$, and it is clear that $U \cap V = \emptyset$. Hence, $H_1$ $\&$ $H_2$ are separated by disjoint $SC^*$-open sets, which establishes that $X$ is $SC^*$-normal.
\end{flushleft}

\section{Functions Associated with {$SC^*$}-Normal Spaces}

\subsection{Definition} A function $f : X \rightarrow Y $ is called.\\
$(1)$ \textbf{R-map}\cite{carnahan1974some} if $f^{-1}(V)$ is regular open in $X$ for every regular open set $V$ of $Y$.\\ 
$(2)$ \textbf{completely continuous}\cite{arya1974strongly} if $f^{-1}(V)$ is regular open in $X$ for every open set $V$ of $Y$.\\
$(3)$ \textbf{rc-continuous}\cite{jankovic1985note} if for each regular closed set $F$ in $Y, f^{-1}(F)$ is regular closed in $X$.

\subsection{Definition} A function $f:X \rightarrow Y$ is called.\\
$(1)$ \textbf{strongly $SC^*$-open} if $f(U)\in SC^*O(Y)$ $\forall$ $U\in SC^*O(X)$.\\
$(2)$ \textbf{strongly $SC^*$-closed} if $f(U)\in SC^*C(Y)$ $\forall$ $U\in SC^*C(X)$.\\
$(3)$ \textbf{almost $SC^*$-irresolute} if $\forall$ $x$ in $X$ and each $SC^*$-nbd $V$ of $f(x)$,\\ $SC^*$-$cl(f^{-1}(V))$ is a $SC^*$-nbd of $x$.

\subsection{Theorem.}
A function \(X \) to \(Y\) is strongly \( SC^* \)-closed iff for every subset \( A \subseteq Y \) and every \( SC^* \)-open set \( U \subseteq X \) with \( f^{-1}(A) \subseteq U \), $\exists$  \( SC^* \)-open set \( V \subseteq Y \) containing \( A \) s.t., 
\[
f^{-1}(V) \subseteq U.
\]

\begin{flushleft}
\textbf{Proof.}
\textbf{($\Rightarrow$)} Assume that \(X \) to \(Y\) is a strongly \( SC^* \)-closed function. Let \( A \subseteq Y \), suppose \( U \in SC^*O(X) \) is an \( SC^* \)-open set s.t., \( f^{-1}(A) \subseteq U \). Define
\[
V = Y \setminus f(X \setminus U).
\]
Then \( V \) is \( SC^* \)-open in \( Y \) because \( X \setminus U \) is \( SC^* \)-closed in \( X \), $\&$ the image under a strongly \( SC^* \)-closed function is also \( SC^* \)-closed in \( Y \). Thus, its complement \( V \) is \( SC^* \)-open. Clearly, \( A \subseteq V \) $\&$ \( f^{-1}(V) \subseteq U \), as required.

\vspace{2mm}

\textbf{($\Leftarrow$)} Conversely, suppose that for every subset \( A \subseteq Y \) $\&$ each \( SC^* \)-open set \( U \subseteq X \) with \( f^{-1}(A) \subseteq U \), $\exists$  \( SC^* \)-open set \( V \subseteq Y \) containing \( A \) s.t., \( f^{-1}(V) \subseteq U \).

To prove that \( f \) is strongly \( SC^* \)-closed, let \( K \subseteq X \) be an \( SC^* \)-closed set. Then \( X \setminus K \in SC^*O(X) \), and since
\[
f^{-1}(Y \setminus f(K)) \subseteq X \setminus K,
\]
by assumption, $\exists$ a \( SC^* \)-open set \( V \subseteq Y \) s.t., \( Y \setminus f(K) \subseteq V \) $\&$ \( f^{-1}(V) \subseteq X \setminus K \). It follows that
\[
f(K) \supseteq Y \setminus V \quad \text{and} \quad K \subseteq f^{-1}(Y \setminus V).
\]
Thus, \( f(K) = Y \setminus V \), which is \( SC^* \)-closed in \( Y \) since \( V \) is \( SC^* \)-open. Therefore, \( f \) is strongly \( SC^* \)-closed.
\end{flushleft}

\subsection{Lemma.}
Let \(X \) to \(Y\) be a function between topological spaces. The following statements are equivalent:

\begin{enumerate}
    \item[(1)] The function \( f \) is almost \( SC^* \)-irresolute; that is, the preimage of every \( SC^* \)-open set in \( Y \) is \( SC^* \)-locally contained in an \( SC^* \)-open set in \( X \).
    
    \item[(2)] For every \( SC^* \)-open set \( V \subseteq Y \), the inclusion  holds
    \[
    f^{-1}(V) \subseteq SC^*\text{-}int\big(SC^*\text{-}cl(f^{-1}(V))\big)
    \].
\end{enumerate}

\subsection{Theorem.} 
A function \(X \) to \(Y\) is almost \( SC^* \)-irresolute iff 
\[
SC^*\text{-}cl(U) \subseteq SC^*\text{-}cl(f(U))
\]
for every \( U \in SC^*O(X) \).

\begin{flushleft}
\textbf{Proof.}
\textbf{($\Rightarrow$)} Suppose \( f \) is almost \( SC^* \)-irresolute. Let \( U \in SC^*O(X) \), $\&$ assume \( y \notin SC^*\text{-}cl(f(U)) \). Then $\exists$ a \( SC^* \)-open set \( V \subseteq Y \) s.t. \( y \in V \) $\&$ \( V \cap f(U) = \emptyset \). It follows that \( f^{-1}(V) \cap U = \emptyset \). 

Since \( U \in SC^*O(X) \), and by the assumption that \( f \) is almost \( SC^* \)-irresolute, we have
\[
f^{-1}(V) \subseteq SC^*\text{-}int(SC^*\text{-}cl(f^{-1}(V))),
\]
which implies 
\[
SC^*\text{-}int(SC^*\text{-}cl(f^{-1}(V))) \cap U = \emptyset.
\]
Thus, taking closures, we get
\[
SC^*\text{-}int(SC^*\text{-}cl(f^{-1}(V))) \cap SC^*\text{-}cl(U) = \emptyset,
\]
which implies 
\[
f^{-1}(V) \cap SC^*\text{-}cl(U) = \emptyset.
\]
Therefore, \( V \cap f(SC^*\text{-}cl(U)) = \emptyset \), $\&$ so \( y \notin f(SC^*\text{-}cl(U)) \). This contradiction shows that \( SC^*\text{-}cl(U) \subseteq SC^*\text{-}cl(f(U)) \).

\vspace{2mm}

\textbf{($\Leftarrow$)} Conversely, assume that for every \( U \in SC^*O(X) \), we have 
\[
SC^*\text{-}cl(U) \subseteq SC^*\text{-}cl(f(U)).
\]
Let \( V \in SC^*O(Y) \), and define 
\[
M = X \setminus SC^*\text{-}cl(f^{-1}(V)).
\]
Then \( M \in SC^*O(X) \). By the hypothesis,
\[
SC^*\text{-}cl(M) \subseteq f^{-1}(SC^*\text{-}cl(f(M))) \subseteq f^{-1}(SC^*\text{-}cl(f(X \setminus f^{-1}(V)))).
\]
Since \( f(X \setminus f^{-1}(V)) \subseteq Y \setminus V \), it follows that
\[
f^{-1}(SC^*\text{-}cl(Y \setminus V)) = f^{-1}(Y \setminus V) = X \setminus f^{-1}(V).
\]
Hence,
\[
X \setminus SC^*\text{-}int(SC^*\text{-}cl(f^{-1}(V))) = SC^*\text{-}cl(M) \subseteq X \setminus f^{-1}(V),
\]
implying
\[
f^{-1}(V) \subseteq SC^*\text{-}int(SC^*\text{-}cl(f^{-1}(V))).
\]
Therefore, by the previous \textbf{lemma 3.4}, \( f \) is almost \( SC^* \)-irresolute.
\end{flushleft}

\subsection{Theorem.} 
Let \(X \) to \(Y\) be a strongly \( SC^* \)-open, cont., $\&$ almost-\( SC^* \)-irresolute surjective function, where \( X \) is a \( SC^* \)-normal topological space. Then \( Y \) is also \( SC^* \)-normal.

\begin{flushleft}
\textbf{Proof.} Let \( A \subseteq Y \) be closed $\&$ \( B \subseteq Y \) be open s.t., \( A \subseteq B \). Since \( f \) is cont., the preimage \( f^{-1}(A) \) is closed in \( X \), and \( f^{-1}(B) \) is open. Thus,
\[
f^{-1}(A) \subseteq f^{-1}(B).
\]

As \( X \) is \( SC^* \)-normal, $\exists$ a \( SC^* \)-open set \( U \subseteq X \) s.t.,
\[
f^{-1}(A) \subseteq U \subseteq SC^*\text{-}cl(U) \subseteq f^{-1}(B).
\]

Applying \( f \) to this chain of inclusions, we obtain
\[
f(f^{-1}(A)) \subseteq f(U) \subseteq f(SC^*\text{-}cl(U)).
\]

Since \( f \) is almost \( SC^* \)-irresolute, it follows that
\[
f(SC^*\text{-}cl(U)) \subseteq SC^*\text{-}cl(f(U)).
\]

Also, as \( f \) is surjective, we have \( f(f^{-1}(A)) = A \) $\&$ \( f(f^{-1}(B)) \subseteq B \). Combining these,
\[
A \subseteq f(U) \subseteq SC^*\text{-}cl(f(U)) \subseteq B.
\]

Thus, for any closed set \( A \subseteq Y \) $\&$ any open set \( B \) containing \( A \), we have identified a \( SC^* \)-open set \( f(U) \) in \( Y \) s.t.,
\[
A \subseteq f(U) \subseteq SC^*\text{-}cl(f(U)) \subseteq B,
\]
which proves that by \textbf{Theorem 3.3} \( Y \) is \( SC^* \)-normal.

\end{flushleft}

\subsection{Theorem.} 
Let \(X \) to \(Y\) be a strongly \( SC^* \)-closed, cont., $\&$ surjective function from a \( SC^* \)-normal space \( X \) onto a space \( Y \). Then \( Y \) is also \( SC^* \)-normal.

\begin{flushleft}
\textbf{Proof.} Let \( M_1, M_2 \subseteq Y \) be two disj., closed sets. Since \( f \) is cont., the preimages \( f^{-1}(M_1) \) $\&$ \( f^{-1}(M_2) \) are closed in \( X \), and they remain disj., As \( X \) is \( SC^* \)-normal, $\exists$ disj., \( SC^* \)-open sets \( U \) $\&$ \( V \) in \( X \) s.t,
\[
f^{-1}(M_1) \subseteq U \quad \text{and} \quad f^{-1}(M_2) \subseteq V.
\]

Since \( f \) is strongly \( SC^* \)-closed and surjective, by the corresponding characterization theorem (e.g., \textbf{Theorem 3.3}), $\exists$ \( SC^* \)-open sets \( A, B \subseteq Y \) s.t.,
\[
M_1 \subseteq A, \quad M_2 \subseteq B, \quad f^{-1}(A) \subseteq U, \quad \text{and} \quad f^{-1}(B) \subseteq V.
\]

Given that \( U \cap V = \emptyset \), it follows that \( f^{-1}(A) \cap f^{-1}(B) = \emptyset \), which implies \( A \cap B = \emptyset \) because \( f \) is surjective. Therefore, \( A \) $\&$ \( B \) are disj., \( SC^* \)-open sets in \( Y \) containing \( M_1 \) and \( M_2 \), respectively.

This proves that \( Y \) is \( SC^* \)-normal.

\end{flushleft}

\section{Generalized $SC^*$-closed functions}

\subsection{Definition}
Let \( f: X \rightarrow Y \) be a function between topological spaces. Then:

\begin{enumerate}
    \item \( f \) is said to be \textbf{\( SC^* \)-closed}\cite{chandrakala2024closed} if for every closed subset \( A \subseteq X \), the image \( f(A) \) is \( SC^* \)-closed in \( Y \).
    
    \item \( f \) is called \textbf{\( SC^*g \)-closed} if for every closed subset \( A \subseteq X \), the image \( f(A) \) is \( SC^*g \)-closed in \( Y \).
    
    \item \( f \) is called \textbf{\( gSC^* \)-closed} if for every closed subset \( A \subseteq X \), the image \( f(A) \) is \( gSC^* \)-closed in \( Y \).
\end{enumerate}

\subsection{Definition}

Let \( f: X \rightarrow Y \) be a function between two topological spaces. Then:

\begin{enumerate}
    \item \( f \) is called \textbf{quasi-\( SC^* \)-closed} if for every \( SC^* \)-closed set \( A \subseteq X \), the image \( f(A) \) is a closed subset of \( Y \).
    
    \item \( f \) is said to be \textbf{\( SC^* \)-\( SC^*g \)-closed} if \( f(A) \) is \( SC^*g \)-closed in \( Y \) for each \( SC^* \)-closed set \( A \subseteq X \).
    
    \item \( f \) is called \textbf{\( SC^* \)-\( gSC^* \)-closed} if for every \( SC^* \)-closed set \( A \subseteq X \), the image \( f(A) \) is \( gSC^* \)-closed in \( Y \).
    
    \item \( f \) is referred to as \textbf{almost-\( gSC^* \)-closed} if the image of every regular closed set in \( X \) under \( f \) is \( gSC^* \)-closed in \( Y \).
\end{enumerate}

\subsection{Definition}
A function \( f: X \rightarrow Y \) is said to be \textbf{\( SC^* \)-\( gSC^* \)-continuous} if the preimage of every \( SC^* \)-closed set in \( Y \) is a \( gSC^* \)-closed set in \( X \).

\subsection{Definition}
A function \( f: X \rightarrow Y \) is said to be \textbf{\( SC^* \)-irresolute} if the preimage of every \( SC^* \)-open set in \( Y \) is also an \( SC^* \)-open set in \( X \); that is, 
\[
f^{-1}(V) \in SC^*O(X) \quad \text{for all} \quad V \in SC^*O(Y).
\]

\subsection{Theorem}
Let \(X \) to \(Y\) $\&$ \(Y \) to \(Z\) be two functions. Then the composition \( g \circ f: X \rightarrow Z \) satisfies the following properties:

\begin{enumerate}
    \item The composition \( g \circ f \) is \( SC^*\)-\(gSC^* \)-closed if \( f \) is \( SC^*\)-\(gSC^* \)-closed $\&$ \( g \) is both cont., $\&$ \( SC^*\)-\(gSC^* \)-closed.
    
    \item The composition \( g \circ f \) is \( SC^*\)-\(gSC^* \)-closed if \( f \) is strongly \( SC^* \)-closed $\&$ \( g \) is \( SC^*\)-\(gSC^* \)-closed.
    
    \item The composition \( g \circ f \) is \( SC^*\)-\(gSC^* \)-closed if \( f \) is quasi \( SC^* \)-closed $\&$ \( g \) is \( gSC^* \)-closed.
\end{enumerate}

\subsection{Theorem}

Let \(X \) to \(Y\) $\&$ \(Y \) to \(Z\) be two functions. If the composition \( g \circ f: X \rightarrow Z \) is\\ \( SC^*\)-\(gSC^* \)-closed and \( f \) is a surjective and \( SC^* \)-irresolute function, then \( g \) is \( SC^*\)-\(gSC^* \)-closed.

\begin{flushleft}
\textbf{Proof.} Consider any \( SC^* \)-closed set \( K \subseteq Y \). Since \( f \) is \( SC^* \)-irresolute and surjective, the preimage \( f^{-1}(K) \) is \( SC^* \)-closed in \( X \). Given that \( g \circ f \) is \( SC^*\)-\(gSC^* \)-closed, the image \( (g \circ f)(f^{-1}(K)) = g(K) \) must be \( gSC^* \)-closed in \( Z \). Thus, \( g \) maps \( SC^* \)-closed sets in \( Y \) to \( gSC^* \)-closed sets in \( Z \), and so \( g \) is \( SC^*\)-\(gSC^* \)-closed.
\end{flushleft}

\subsection{Lemma.} A mapping $X$ to \(Y\) is $SC^*$-$gSC^*$-closed iff for each subset $B\subset Y$ $\&$ each\\ $U\in SC^*O(X)$ containing $f^{-1}(B)$, $\exists$ a $gSC^*$-open set of $Y$ s.t., $B\subset V$ $\&$ $f^{-1}(V)\subset U$.

\textbf{Proof.}
$(\Rightarrow)$: Assume that \( f \) is a \( SC^* \)-\( gSC^* \)-closed function. Let \( B \subseteq Y \), $\&$ suppose\\ \( f^{-1}(B) \subseteq U \), where \( U \in SC^*O(X) \). Define the set \( V = Y \setminus f(X \setminus U) \). Since \( X \setminus U \) is \( SC^* \)-closed $\&$ \( f \) maps \( SC^* \)-closed sets to \( gSC^* \)-closed sets, the image \( f(X \setminus U) \) is \( gSC^* \)-closed in \( Y \). Therefore, \( V \) is \( gSC^* \)-open in \( Y \), and it satisfies \( B \subseteq V \) $\&$ \( f^{-1}(V) \subseteq U \), as required.

\vspace{2mm}

$(\Leftarrow)$: Let \( K \) be a \( SC^* \)-closed subset of \( X \). Then the complement \( X \setminus K \) is \( SC^* \)-open, and we observe that \( f^{-1}(Y \setminus f(K)) \subseteq X \setminus K \). By the given condition,  $\exists$ a \( gSC^* \)-open set \( V \subseteq Y \) s.t., \( Y \setminus f(K) \subseteq V \) $\&$ \( f^{-1}(V) \subseteq X \setminus K \). It follows that \( f(K) \supseteq Y \setminus V \), $\&$ since \( Y \setminus V \) is \( gSC^* \)-closed, we conclude that \( f(K) \) is also \( gSC^* \)-closed. Hence, \( f \) is \( SC^* \)-\( gSC^* \)-closed.

\subsection{Theorem.} 
Let \(X \) to \(Y\) be a cont., function that is also \( SC^* \)-\( gSC^* \)-closed. Then the image of every \( gSC^* \)-closed set in \( X \) under \( f \) is \( gSC^* \)-closed in \( Y \).

\begin{flushleft}
\textbf{Proof.} 
Consider a set \( H \subseteq X \) that is \( gSC^* \)-closed, $\&$ let \( V \subseteq Y \) be an open set s.t, \( f(H) \subseteq V \). By the continuity of \( f \), the preimage \( f^{-1}(V) \) is open in \( X \) and contains \( H \).
Since \( H \) is \( gSC^* \)-closed, we know that \( SC^*\text{-}cl(H) \subseteq f^{-1}(V) \). Applying the function \( f \) gives:
\[
f(SC^*\text{-}cl(H)) \subseteq f(f^{-1}(V)) \subseteq V.
\]

Now, because \( SC^*\text{-}cl(H) \) is a \( SC^* \)-closed set $\&$ \( f \) is \( SC^* \)-\( gSC^* \)-closed, it follows that \( f(SC^*\text{-}cl(H)) \) is \( gSC^* \)-closed in \( Y \). Therefore:
\[
SC^*\text{-}cl(f(H)) \subseteq SC^*\text{-}cl(f(SC^*\text{-}cl(H))) \subseteq V.
\]

Since this holds for every open set \( V \) containing \( f(H) \), it follows that \( f(H) \) is \( gSC^* \)-closed in \( Y \).
\end{flushleft}

\subsection{Remark.}  
Every $SC^*$-irresolute function is $SC^*$-$gSC^*$-cont., but the converse does not necessarily hold.

\subsection{Theorem.}  
A function \(X \) to \(Y\) is $SC^*$-$gSC^*$-cont., iff for every $V \in SC^*O(Y)$, the preimage $f^{-1}(V)$ is $gSC^*$-open in $X$.

\subsection{Theorem.}  
If a function \(X \) to \(Y\) is $SC^*$-$gSC^*$-cont., then the preimage of every $gSC^*$-closed set in $Y$ is $gSC^*$-closed in $X$.

\begin{flushleft}
\textbf{Proof.}  
Let $K$ be a $gSC^*$-closed subset of $Y$, $\&$ let $U$ be an open set in $X$ s.t., $f^{-1}(K) \subset U$. Define $V = Y \setminus f(X \setminus U)$. Then $V$ is open in $Y$ $\&$ satisfies $K \subset V$ $\&$ $f^{-1}(V) \subset U$.  
Since $SC^*$-$cl(K) \subset V$, it follows that:
\[
f^{-1}(K) \subset f^{-1}(SC^*\text{-}cl(K)) \subset f^{-1}(V) \subset U.
\]
Because $f$ is $SC^*$-$gSC^*$-cont., $f^{-1}(SC^*$-$cl(K))$ is $gSC^*$-closed in $X$. Therefore,
\[
SC^*\text{-}cl(f^{-1}(K)) \subset SC^*\text{-}cl(f^{-1}(SC^*\text{-}cl(K))) \subset U,
\]
which implies that $f^{-1}(K)$ is $gSC^*$-closed in $X$.
\end{flushleft}

\subsection{Corollary.}  
If \(X \) to \(Y\) is a closed $\&$ $SC^*$-irresolute function, then for every $gSC^*$-closed set $K\in Y$, the preimage $f^{-1}(K)$ is $gSC^*$-closed in $X$.

\subsection{Theorem.}  
Let \(X \) to \(Y\) be a bijective, open $\&$ $SC^*$-$gSC^*$-cont., function. Then for any $gSC^*$-closed set $K\in Y$, the preimage $f^{-1}(K)$ is $gSC^*$-closed in $X$.

\begin{flushleft}
\textbf{Proof.}  
Assume $K$ is a $gSC^*$-closed subset of $Y$, and let $U$ be an open subset of $X$ s.t. $f^{-1}(K) \subset U$. Since $f$ is a surjective open map, we have:
\[
K = f(f^{-1}(K)) \subset f(U),
\]
where $f(U)$ is open in $Y$. This implies that:
\[
SC^*\text{-}cl(K) \subset f(U).
\]
Due to injectivity of $f$, we get:
\[
f^{-1}(K) \subset f^{-1}(SC^*\text{-}cl(K)) \subset f^{-1}(f(U)) = U.
\]
Now, because $f$ is $SC^*$-$gSC^*$-cont., the preimage $f^{-1}(SC^*\text{-}cl(K))$ is $gSC^*$-closed in $X$. Consequently,
\[
SC^*\text{-}cl(f^{-1}(K)) \subset SC^*\text{-}cl(f^{-1}(SC^*\text{-}cl(K))) \subset U,
\]
which confirms that $f^{-1}(K)$ is $gSC^*$-closed in $X$.
\end{flushleft}

\subsection{Theorem.}  
Let \(X \) to \(Y\) $\&$ \(Y \) to \(Z\) be two functions. If the composition $g \circ f: X \rightarrow Z$ is $SC^*$-$gSC^*$-closed $\&$ $g$ is a bijective, open, and $SC^*$-$gSC^*$-cont., function, then $f$ is also $SC^*$-$gSC^*$-closed.

\begin{flushleft}
\textbf{Proof.}  
Consider any $SC^*$-closed set $H$ in $X$. Since $g \circ f$ is $SC^*$-$gSC^*$-closed by assumption, the image $(g \circ f)(H)$ is $gSC^*$-closed in $Z$. Because $g$ is bijective, we can write:
\[
f(H) = g^{-1}((g \circ f)(H)).
\]
Since $g$ is $SC^*$-$gSC^*$-cont., $\&$ open, the preimage of a $gSC^*$-closed set under $g$ remains $gSC^*$-closed in $Y$. Therefore, $f(H)$ is $gSC^*$-closed in $Y$. This confirms that $f$ is $SC^*$-$gSC^*$-closed.
\end{flushleft}

\subsection{Theorem.} 
Let \(X \) to \(Y\) $\&$ \(Y \) to \(Z\) be two functions s.t. the composition $g \circ f: X \rightarrow Z$ is $SC^*$-$gSC^*$-closed. If $g$ is an injective, $\&$ $SC^*$-$gSC^*$-cont., function, so $f$ is also $SC^*$-$gSC^*$-closed.

\begin{flushleft}
\textbf{Proof.}  
Consider a set $H \in SC^*C(X)$, meaning that $H$ is $SC^*$-closed in $X$. Since the composition $g \circ f$ is $SC^*$-$gSC^*$-closed, the image $(g \circ f)(H)$ is $gSC^*$-closed in $Z$.

Because $g$ is injective, we have:
\[
f(H) = g^{-1}((g \circ f)(H)).
\]

Now, using the fact that $g$ is both closed $\&$ $SC^*$-$gSC^*$-cont., the preimage of a $gSC^*$-closed set under $g$ remains $gSC^*$-closed in $Y$. Therefore, $f(H)$ is $gSC^*$-closed in $Y$, which confirms that $f$ is $SC^*$-$gSC^*$-closed.
\end{flushleft}

\section{Preservation theorems and other characterizations of $SC^*$-normal spaces}
\subsection{Theorem.} 
For a topological space $X$, the following conditions are equivalent:

\begin{enumerate}[(a)]
    \item $X$ is $SC^*$-normal.
    \item Any two disjoint closed subsets $A$ $\&$ $B$ of $X$ can be separated by disjoint $gSC^*$-open sets $U$ $\&$ $V$ s.t., $A \subset U$ $\&$ $B \subset V$.
    \item Given any closed set $A$ $\&$ open set $B$ with $A \subset B$,  $\exists$ a $gSC^*$-open set $U$ s.t.,\\ $\text{cl}(A) \subset U \subset SC^*\text{-cl}(U) \subset B$.
    \item For every closed set $A$ $\&$ each $g$-open set $B$ containing $A$, $\exists$ an $SC^*$-open set $U$ s.t. \\$A \subset U \subset SC^*\text{-cl}(U) \subset \text{int}(B)$.
    \item For any closed set $A$ $\&$ $g$-open set $B$ containing $A$,  $\exists$ a $gSC^*$-open set $G$ s.t.\\ $A \subset G \subset SC^*\text{-cl}(G) \subset \text{int}(B)$.
    \item If $A$ is $g$-closed $\&$ $B$ is an open set with $A \subset B$, then $\exists$ an $SC^*$-open set $U$ s.t.\\ $\text{cl}(A) \subset U \subset SC^*\text{-cl}(U) \subset B$.
    \item If $A$ is $g$-closed $\&$ $B$ is open with $A \subset B$, then  $\exists$ a $gSC^*$-open set $G$ s.t.,\\ $\text{cl}(A) \subset G \subset SC^*\text{-cl}(G) \subset B$.
\end{enumerate}

\subsection*{Proof.} 

\noindent
\textbf{(a) $\Leftrightarrow$ (b) $\Leftrightarrow$ (c):}  
These equivalences follow naturally since any $SC^*$-open set is also $gSC^*$-open, so the separation and approximation conditions preserve the same structure.

\noindent
\textbf{(d) $\Rightarrow$ (e) $\Rightarrow$ (c) $\&$ (f) $\Rightarrow$ (g) $\Rightarrow$ (c):}  
These implications hold because every closed (resp. open) set is also $g$-closed (resp. $g$-open), making the conditions more general.

\noindent
\textbf{(c) $\Rightarrow$ (e):}  
Let $A$ be closed $\&$ $B$ a $g$-open set with $A \subset B$. Since $B$ is $g$-open, it follows that $A \subset \text{int}(B)$. Then by \textbf{(c)}, $\exists$ a $gSC^*$-open set $U$ s.t, $A \subset U \subset SC^*\text{-cl}(U) \subset \text{int}(B)$.

\noindent
\textbf{(e) $\Rightarrow$ (d):}  
Given $A$ closed $\&$ $B$ a $g$-open nbd of $A$, let $G$ be a $gSC^*$-open set satisfying\\ $A \subset G \subset SC^*\text{-cl}(G) \subset \text{int}(B)$. Since $G$ is $gSC^*$-open, its $SC^*$-interior $U = SC^*\text{-int}(G)$ is an $SC^*$-open set with $A \subset U \subset SC^*\text{-cl}(U) \subset \text{int}(B)$.

\noindent
\textbf{(c) $\Rightarrow$ (g):}  
If $A$ is $g$-closed $\&$ $B$ is open with $A \subset B$, then $\text{cl}(A) \subset B$. From \textbf{(c)}, $\exists$ a $gSC^*$-open set $U$ with $\text{cl}(A) \subset U \subset SC^*\text{-cl}(U) \subset B$, satisfying \textbf{(g)}.

\noindent
\textbf{(g) $\Rightarrow$ (f):}  
Let $A$ be $g$-closed $\&$ $B$ open with $A \subset B$. Then a $gSC^*$-open set $G$ $\exists$ s.t.,\\ $\text{cl}(A) \subset G \subset SC^*\text{-cl}(G) \subset B$. As $G$ is $gSC^*$-open and $\text{cl}(A) \subset G$, we obtain $\text{cl}(A) \subset SC^*\text{-int}(G)$. Setting $U = SC^*\text{-int}(G)$, we get an $SC^*$-open set satisfying $\text{cl}(A) \subset U \subset SC^*\text{-cl}(U) \subset B$.

\subsection{Theorem.} 
Let \(X \) to \(Y\) be a cont., quasi $SC^*$-closed surjective function. If $X$ is $SC^*$-normal, then $Y$ is a normal space.

\begin{flushleft}
\textbf{Proof.} 
Assume $M_1$ $\&$ $M_2$ are two disj., closed subsets of $Y$. Since $f$ is cont., their preimages $f^{-1}(M_1)$ $\&$ $f^{-1}(M_2)$ are closed and disj., in $X$. Given that $X$ is $SC^*$-normal, $\exists$ disj., $SC^*$-open sets $U_1$ $\&$ $U_2$ in $X$ s.t., $f^{-1}(M_1) \subset U_1$ $\&$ $f^{-1}(M_2) \subset U_2$.
Now, define $V_i = Y \setminus f(X \setminus U_i)$ for $i = 1, 2$. Because $f$ is quasi $SC^*$-closed, the sets $V_1$ $\&$ $V_2$ are open in $Y$. Moreover, we have $M_i \subset V_i$ $\&$ $f^{-1}(V_i) \subset U_i$ for each $i$. Since $U_1 \cap U_2 = \emptyset$, $\&$ $f$ is surjective, it follows that $V_1 \cap V_2 = \emptyset$. 
Hence, the disj., closed sets $M_1$ $\&$ $M_2$ in $Y$ can be separated by disj., open sets, which implies that $Y$ is normal.
\end{flushleft}

\subsection{Lemma.} 
A subset \( A \) of a topological space \( X \) is said to be \( gSC^* \)-open iff, for every closed set \( F \subseteq A \), it follows that \( F \subseteq SC^*\text{-}int(A) \).

\subsection{Theorem.} 
Let \(X \) to \(Y\) be a closed, injective function that is also \( SC^*\text{-}gSC^* \)-cont.. If the codomain \( Y \) is \( SC^* \)-normal, then the domain \( X \) must also be \( SC^* \)-normal.

\begin{flushleft}
\textbf{Proof.} 
Consider two disj., closed subsets \( N_1 \) $\&$ \( N_2 \) of \( X \). Since \( f \) is a closed injection, the images \( f(N_1) \) $\&$ \( f(N_2) \) are closed and non-overlapping in \( Y \). Because \( Y \) is \( SC^* \)-normal, there exist disj., \( SC^* \)-open sets \( V_1 \) $\&$ \( V_2 \) in \( Y \) s.t., \( f(N_1) \subset V_1 \) $\&$ \( f(N_2) \subset V_2 \).
Due to the \( SC^*\text{-}gSC^* \)-continuity of \( f \), the inverse images \( f^{-1}(V_1) \) $\&$ \( f^{-1}(V_2) \) are \( gSC^* \)-open in \( X \) and contain \( N_1 \) $\&$ \( N_2 \), respectively. Define \( U_i = SC^*\text{-}int(f^{-1}(V_i)) \) for \( i = 1, 2 \). Then \( U_1 \) $\&$ \( U_2 \) are \( SC^* \)-open, disj., $\&$ each \( N_i \subset U_i \), establishing that \( X \) is \( SC^* \)-normal.
\end{flushleft}

\subsection{Corollary.} 
If \(X \) to \(Y\) is a closed, \( SC^* \)-irresolute, $\&$ injective function, and if \( Y \) is \( SC^* \)-normal, then \( X \) must also be \( SC^* \)-normal.

\begin{flushleft}
\textbf{Proof.} 
This result follows directly since every \( SC^* \)-irresolute function is inherently \( SC^*\text{-}gSC^* \)-continuous.
\end{flushleft}

\subsection{Lemma.} 
A mapping \(X \) to \(Y\) is almost-\( gSC^* \)-closed iff for every subset \( B \subseteq Y \) $\&$ every \(r\)-open set \( U \in RO(X) \) satisfying \( f^{-1}(B) \subseteq U \), $\exists$ a \( gSC^* \)-open set \( V \subseteq Y \) s.t., \( B \subseteq V \) $\&$ \( f^{-1}(V) \subseteq U \).

\subsection{Lemma.} 
Let \(X \) to \(Y\) be an almost \( gSC^* \)-closed function. Then for any closed subset \( M \subseteq Y \) $\&$ any \(r\)-open set \( U \in RO(X) \) with \( f^{-1}(M) \subseteq U \), $\exists$ an \( SC^* \)-open set \( V \subseteq Y \) s.t., \( M \subseteq V \) $\&$ \( f^{-1}(V) \subseteq U \).

\subsection{Theorem.} 
Let \(X \) to \(Y\) be a cont., almost \( gSC^* \)-closed, $\&$ surjective function. If the space \( X \) is normal, then the space \( Y \) is \( SC^* \)-normal.

\begin{flushleft}
\textbf{Proof.} 
Let \( M_1 \) $\&$ \( M_2 \) be two disj., closed subsets of \( Y \). By continuity of \( f \), their preimages \( f^{-1}(M_1) \) $\&$ \( f^{-1}(M_2) \) are disj., closed subsets in \( X \). Since \( X \) is normal, $\exists$ disj., open sets \( U_1, U_2 \subset X \) s.t., \( f^{-1}(M_i) \subseteq U_i \) for each \( i = 1, 2 \).
Define \( G_i = \text{int}(\text{cl}(U_i)) \), which belongs to the regular open sets of \( X \), denoted by \( RO(X) \), $\&$ satisfies \( f^{-1}(M_i) \subseteq G_i \) with \( G_1 \cap G_2 = \emptyset \). 
Since \( f \) is almost \( gSC^* \)-closed, for each \( i \), $\exists$ a \( gSC^* \)-open set \( V_i \subset Y \) s.t., \( M_i \subseteq V_i \) $\&$ \( f^{-1}(V_i) \subseteq G_i \). 
Because \( f \) is surjective $\&$ \( G_1 \cap G_2 = \emptyset \), it follows that \( V_1 \cap V_2 = \emptyset \). Thus, \( M_1 \) $\&$ \( M_2 \) are contained in disj., \( SC^* \)-open sets in \( Y \), proving that \( Y \) is \( SC^* \)-normal.
\end{flushleft}

\subsection{Corollary.} 
Let \(X \) to \(Y\) be a cont. $\&$ \( SC^* \)-closed surjective function. If \( X \) is a normal topological space, then \( Y \) is \( SC^* \)-normal.

\begin{flushleft}
\textbf{Proof.} 
Let \( M_1 \) $\&$ \( M_2 \) be two disj., closed subsets of \( Y \). Since \( f \) is cont., the preimages \( f^{-1}(M_1) \) $\&$ \( f^{-1}(M_2) \) are closed in \( X \) $\&$ also disj., Because \( X \) is normal,  $\exists$ disj., open sets \( U_1, U_2 \subset X \) s.t., \( f^{-1}(M_i) \subset U_i \) for \( i = 1, 2 \).
Define \( G_i = \operatorname{int}(\overline{U_i}) \), which belongs to the regular open subsets of \( X \), and note that \( f^{-1}(M_i) \subset G_i \) $\&$ \( G_1 \cap G_2 = \varnothing \). Since \( f \) is \( SC^* \)-closed, it is also almost-\( gSC^* \)-closed.
By the corresponding lemma for almost \( gSC^* \)-closed functions, for each \( i = 1, 2 \), $\exists$ a \( SC^* \)-open set \( V_i \subset Y \) s.t., \( M_i \subset V_i \) $\&$ \( f^{-1}(V_i) \subset G_i \). Since \( G_1 \cap G_2 = \varnothing \), it follows that \( V_1 \cap V_2 = \varnothing \).
Thus, we have disj., \( SC^* \)-open sets \( V_1 \) $\&$ \( V_2 \) in \( Y \) containing \( M_1 \) $\&$ \( M_2 \), respectively. Therefore, \( Y \) is \( SC^* \)-normal.
\end{flushleft}


\bibliographystyle{amsplain}
\bibliography{references}

\end{document}